\documentclass{article}

\usepackage{amsmath}
\usepackage{amsthm}
\usepackage{amssymb}
\usepackage{amscd}
\usepackage{ifthen}
\usepackage{epsfig,graphics,graphicx}
\usepackage{array,tabularx}
\usepackage{url}
\usepackage[usenames]{color}
\usepackage[pdftex,colorlinks=true,urlcolor=blue,linkcolor=blue,plainpages=false,pdfpagelabels]{hyperref}

\newtheorem{theorem}{Theorem}[section]

\newtheorem{definition}[theorem]{Definition}
\newtheorem{proposition}[theorem]{Proposition}
\newtheorem{corollary}[theorem]{Corollary}
\newtheorem{lemma}[theorem]{Lemma}
\newtheorem{fact}[theorem]{Remark}
\newtheorem{exemplu}[theorem]{Example}
\newtheorem{exercise}{Exercise}
\newtheorem{notation}[theorem]{Notation}

\newcommand{\bdfn}{\begin{definition}}
\newcommand{\edfn}{\end{definition}}
\newcommand{\bthm}{\begin{theorem}}
\newcommand{\ethm}{\end{theorem}}
\newcommand{\bprop}{\begin{proposition}}
\newcommand{\eprop}{\end{proposition}}
\newcommand{\bcor}{\begin{corollary}}
\newcommand{\ecor}{\end{corollary}}
\newcommand{\blem}{\begin{lemma}}
\newcommand{\elem}{\end{lemma}}
\newcommand{\bfact}{\begin{fact}}
\newcommand{\efact}{\end{fact}}
\newcommand{\bex}{\begin{exemplu}\begin{rm}}
\newcommand{\eex}{\end{rm}\end{exemplu}}
\newcommand{\bxc}{\begin{exercise}}
\newcommand{\exc}{\end{exercise}}
\newcommand{\bntn}{\begin{notation}}
\newcommand{\entn}{\end{notation}}

\newcommand{\be}{\begin{enumerate}}
\newcommand{\ee}{\end{enumerate}}
\newcommand{\bce}{\begin{center}}
\newcommand{\ece}{\end{center}}
\newcommand{\bi}{\begin{itemize}}
\newcommand{\ei}{\end{itemize}}
\newcommand{\bt}{\begin{tabular}}
\newcommand{\et}{\end{tabular}}
\newcommand{\beq}{\begin{equation}}
\newcommand{\eeq}{\end{equation}}
\newcommand{\ba}{\begin{array}} 
\newcommand{\ea}{\end{array}}
\newcommand {\bea} {\begin{eqnarray}}
\newcommand {\eea} {\end {eqnarray}}
\newcommand {\bua} {\begin{eqnarray*}}
\newcommand {\eua} {\end {eqnarray*}}

\newcommand{\se}{\subseteq}

\newcommand{\ds}{\displaystyle}

\def\R{{\mathbb R}}
\def\N{{\mathbb N}}
\def\Z{{\mathbb Z}}

\def\R{{\mathbb R}}

\newcommand{\eps}{\varepsilon}

\newcommand{\limn}{\ds\lim_{n\to\infty}}

\newcounter{ct}

\newcommand{\solution}[1]{\begin{proof}#1\end{proof}}

\newcommand{\tnN}{T_{n,N}}

\begin{document}

\title{Effective results on compositions of nonexpansive mappings}
\author{Lauren\c{t}iu Leu\c{s}tean${}^{1}$,  Adriana Nicolae${}^{2,3}$ \\[0.2cm]
\footnotesize ${}^1$ Simion Stoilow Institute of Mathematics of the Romanian Academy, Research unit 5,\\
\footnotesize P. O. Box 1-764, RO-014700 Bucharest, Romania\\[0.1cm]
\footnotesize ${}^2$ Department of Mathematics, Babe\c{s}-Bolyai University, \\
\footnotesize  Kog\u{a}lniceanu 1, 400084 Cluj-Napoca, Romania\\[0.1cm]
\footnotesize ${}^3$ Simion Stoilow Institute of Mathematics of the Romanian Academy,\\
\footnotesize Research group of the project PD-3-0152,\\
\footnotesize P. O. Box 1-764, RO-014700 Bucharest, Romania\\[0.1cm]
\footnotesize E-mails:  Laurentiu.Leustean@imar.ro, anicolae@math.ubbcluj.ro
}

\maketitle

\begin{abstract}
This paper provides uniform bounds on the  asymptotic regularity for iterations associated to a finite family of nonexpansive mappings. 
We obtain our  quantitative results  in the setting of $(r,\delta)$-convex spaces,  a class of geodesic spaces which generalizes metric spaces with a convex geodesic bicombing. 
\\

\noindent {\em MSC:} 47J25; 47H09;  53C23; 03F10.\\

\noindent {\em Keywords:}  Effective rates of asymptotic regularity;  Proof mining; Families of nonexpansive mappings;  
Halpern iteration; Convex geodesic bicombing.

\end{abstract}

\section{Introduction}

Let $X$ be a  Hilbert space, $C\se X$ a closed convex subset, $T_1,\ldots, T_N:C\to C$ (where $N\in\Z_+$) 
a finite family of nonexpansive mappings and $(\lambda_n)$ a sequence in $[0,1]$. 
Given $u\in C$, one can define an iteration starting from $u$ by
\beq
x_0=u, \quad x_{n+1} = \lambda_{n+1}u + (1-\lambda_{n+1})T_{n+1}x_n, \label{Halpern-N} 
\eeq
where $T_n = T_{n \text{ mod } N}$ and the $\text{mod } N$ function takes values in $1, \ldots ,N$. For the special case $N=1$, this iteration coincides 
with the well-known Halpern iteration \cite{Hal67}, whose strong convergence was proved 
by Wittmann \cite{Wit92} under suitable conditions on  $(\lambda_n)$, which are satisfied 
by the natural candidate $\lambda_n=\frac1{n+1}$.

The general iteration  defined by \eqref{Halpern-N} was first studied in Hilbert spaces by 
Lions \cite{Lio77}, who assumed different hypotheses on $(\lambda_n)$, with the drawback 
that $\lambda_n=\frac1{n+1}$ does not satisfy them. Bauschke \cite{Bau96} proved that the iteration 
$(x_n)$ given by \eqref{Halpern-N} converges strongly to the common fixed point of the 
mappings $T_1, \ldots, T_N$ which is nearest to $u$, under the assumptions that 
$\ds F:=\bigcap_{i=1}^N\text{Fix}(T_i)$ is nonempty, $F=\text{Fix}(T_NT_{N-1}\cdots T_1)=\ldots = \text{Fix}(T_1T_N\cdots T_2)=
\text{Fix}(T_{N-1}\cdots T_1T_N)$ and that $(\lambda_n)$ satisfies 
\beq
\ds \lim_{n\to\infty} \lambda_n =0, \quad  \sum_{n=1}^\infty|\lambda_{n+N}-\lambda_n| <\infty \quad \text {and} \quad 
\sum_{n=1}^\infty \lambda_n=\infty. \label{hyp-lambdan}
\eeq

Bauschke's result is a generalization of Wittmanns's theorem to a finite family of mappings, since for $N=1$
conditions \eqref{hyp-lambdan} coincide with the ones used by Wittmann.

Another iteration that will be considered in this paper is the following one associated to any 
nonexpansive mapping $T:C\to C$,
\beq
x_{n+1} = T(\lambda_{n+1}u + (1-\lambda_{n+1})x_n). \label{alternative-iteration}
\eeq
This iteration, studied by Xu \cite{Xu10}, is a discrete version of the approximating curve 
$z_t = T(tu + (1-t)z_t), \, t\in (0,1)$,
analyzed by Combettes and Hirstoaga \cite{ComHir06}. Strong convergence of the iteration 
\eqref{alternative-iteration} was 
established by Xu in the setting of uniformly smooth Banach spaces under appropriate assumptions 
on $(\lambda_n)$, including those given by \eqref{hyp-lambdan}. The strong  convergence 
results of Xu  and Combettes and Hirstoaga were extended in \cite{ColLeuLopMar11} 
to more general approximating curves and iterations. As above, one can define for the iteration \eqref{alternative-iteration}
a cyclic algorithm associated to the finite family of nonexpansive mappings $T_1,\ldots, T_N:C\to C$,
\beq
x_0=u, \quad x_{n+1} = T_{n+1}(\lambda_{n+1}u + (1-\lambda_{n+1})x_n). \label{alternative-iteration-N}
\eeq

A very important concept in the study of the asymptotic behavior of nonlinear iterations is the so-called 
asymptotic regularity, introduced by Browder and Petryshyn \cite{BroPet66} in their study of solutions of 
nonlinear functional equations using Picard iterations: $T:C\to C$ is asymptotically regular if 
$\ds \limn \|T^nx-T^{n+1}x\|=0$ for all $x\in C$. More generally, an 
iteration $(x_n)$ associated to a mapping $T$ is said to be {\em asymptotically regular} if 
$\ds \limn \|x_n-Tx_n\|=0$ for all starting points in $C$.

A natural question is to compute rates of asymptotic regularity for the iteration $(x_n)$, 
i.e. rates of convergence of $(\|x_n-Tx_n\|)$ towards $0$.  For the Halpern iteration this 
was done  in  a series of papers 
\cite{Leu07a,Leu09-Habil,Koh11,KohLeu12a,LeuNic13}, corresponding to different classes of 
spaces.
For the iteration given by  \eqref{alternative-iteration},  such rates were obtained in 
\cite{ColLeuLopMar11}. 

The notion of asymptotic regularity can be extended to sequences $(x_n)$ associated to a 
family of mappings $T_1,\ldots, T_N:C\to C$, as it is the case in this paper. Thus, we say that $(x_n)$ is 
{\em asymptotically regular} if $\ds \limn \|x_n-T_{n+N}\cdots T_{n+1}x_n\|=0$ for all starting points in $C$. 
The following asymptotic regularity result is contained in Bauschke's 
strong convergence proof for the iteration \eqref{Halpern-N}.

\begin{theorem}\label{asymp-reg-N}
Let $X$ be a Hilbert space, $C\se X$  convex,  $T_1,\ldots, T_N:C\to C$  nonexpansive mappings 
and $(\lambda_n)$  a sequence in $[0,1]$ satisfying \eqref{hyp-lambdan}.  Let $(x_n)$ be 
given by \eqref{Halpern-N} and assume that $(x_n)$ is bounded. Then, 
\[\limn  \|x_n-T_{n+N}\cdots T_{n+1}x_n\|=0.\]
\end{theorem}

The main result of this paper is a quantitative version of Theorem \ref{asymp-reg-N}
for both iterations \eqref{Halpern-N} and \eqref{alternative-iteration-N}. In order to get this result we apply methods of proof mining developed by Kohlenbach \cite{Koh08-book} with the aim of obtaining effective and uniform bounds  from proofs where such information is not readily available. 
As a consequence, we provide for the first time effective and uniform rates of asymptotic regularity for the 
iterations \eqref{Halpern-N} and \eqref{alternative-iteration-N}. 

Actually, we obtain our quantitative results in a setting more general than the one of normed 
space. More precisely, we introduce {\em $(r,\delta)$-convex spaces}, a class of metric spaces which also includes Busemann spaces (and, hence, CAT$(0)$ spaces), 
hyperconvex spaces, CAT$(\kappa)$ spaces with $\kappa > 0$, as well as the so-called $W$-hyperbolic spaces (see \cite{Koh05}).
Consequently, even when $N=1$ and so \eqref{Halpern-N} reduces in fact to the Halpern iteration, our results generalize ones obtained previously by 
the authors for CAT$(\kappa)$ spaces with $\kappa>0$ 
\cite{LeuNic13} and by the first author for normed \cite{Leu07a} or $W$-hyperbolic 
spaces \cite{Leu09-Habil}.

\section{$(r,\delta)$-convex spaces}

Let $(X,d)$ be a metric space.  We recall first basic facts in geodesic geometry. Given $x,y\in X$, 
a {\it constant speed geodesic } from $x$ to $y$ is a mapping $\gamma:[0,1]\to X$ such that 
$\gamma(0)=x$, $\gamma(1)=y$ and 
$d(\gamma(s),\gamma(t))=|s-t|d(x,y)$ for all $s,t\in [0,1]$. The image $\gamma([0,1])$ of $\gamma$ is 
a {\it geodesic segment} which joins $x$ and $y$. Note that a geodesic segment from $x$ to $y$ 
is not necessarily unique. Given $r\in(0,\infty]$, we say that $(X,d)$ is a {\it (uniquely) $r$-geodesic space} if every two points 
$x,y \in X$ with $d(x,y)\leq r$ can be joined by a (unique) geodesic segment. For $r=\infty$, 
we say simply that $X$ is a (uniquely) geodesic space.

Let $r\in(0,\infty]$ and $X$ be an $r$-geodesic space. We consider an $r$-geodesic bicombing $\Gamma$ on $X$, that is, a choice of a constant speed geodesic 
$\gamma_{x,y}$ joining $x$ and $y$ for each pair of points $x,y\in X$ with $d(x,y)\leq r$. 
When $r=\infty$, $\Gamma$ is called a geodesic bicombing on $X$. 
If $X$ is uniquely $r$-geodesic, then clearly one can define an  $r$-geodesic bicombing in a unique way. 
For $\gamma_{x,y}\in \Gamma$, we denote by  $[x,y]$  the geodesic segment 
$\gamma_{x,y}([0,1])$. A subset $C$  of $X$ is 
{\it $r$-convex} if $[x,y]\in C$ for all $x,y\in C$ with $d(x,y)\leq r$. Given $t\in [0,1]$, we use the notation $(1-t)x+ty$ for $\gamma_{x,y}(t)$. Then, $d(x,(1-t)x+ty) = td(x,y)$ and 
$d(y,(1-t)x+ty)=(1-t)d(x,y)$. The $r$-geodesic bicombing  is convex if it satisfies
\beq
d((1-t)x+ty, (1-t)x+tz)\leq td(y,z) \label{ineq-Busemann-convex}
\eeq
for all $x,y,z\in X$ with $d(x,y),d(x,z),d(z,y)\leq r$ and all $t\in [0,1]$. 

Normed spaces are obviously geodesic spaces with a convex geodesic bicombing. Another natural 
example are Busemann spaces, which were used for the first 
time by Busemann \cite{Bus55} to give a definition of 
nonpositive curvature in geodesic spaces. Thus, geodesic spaces with the property that 
each point has 
a convex neighborhood which is a Busemann space are `nonpositively curved' spaces in the 
sense of Busemann, who called them $G$-spaces. 
We refer to \cite{Pap05} for a nice exposition of this very important class of geodesic spaces. 
It turns out that Busemann spaces are 
uniquely geodesic spaces with a (unique) convex geodesic bicombing. 

A related example of  metric spaces with a  convex geodesic bicombing are the so-called $W$-hyperbolic spaces, 
defined in \cite{Koh05} as metric spaces together with  a convexity mapping $W:X\times X\times [0,1]\to X$ 
satisfying suitable properties.  As it was remarked in \cite{AriLeuLop12}, Busemann spaces 
are exactly the uniquely geodesic $W$-hyperbolic 
spaces.

It is well-known that a hyperconvex space $X$ also admits a convex geodesic bicombing obtained by embedding $X$ 
isometrically into $\ell^\infty(X)$ and using the existence of a nonexpansive retraction 
from $\ell^\infty(X)$ into $X$ to define the geodesic bicombing through the convex linear 
geodesic bicombing on $\ell^\infty(X)$ (for more details see, for instance, Chapter 13 in \cite{KirSim01}). 

In the following we define a natural generalization of metric spaces with a convex $r$-geodesic bicombing.

\bdfn\label{def-rdelta-Busemann}
Let $r\in(0,\infty]$ and $\delta\in[0,1]$. A metric space $(X,d)$ with an $r$-geodesic bicombing is said to be 
$(r,\delta)$-convex if for all $x,y,z\in X$ with $d(x,y),d(x,z),d(z,y)\leq r$ and all $t\in [0,1]$,
\[
d((1-t)x+ty, (1-t)x+tz)\leq (t+\delta(1-t))d(y,z).
\]
If $r=\infty$, we say that $X$ is $\delta$-convex.
\edfn

An example  of such spaces are CAT$(\kappa)$ spaces with $\kappa > 0$. CAT$(\kappa)$ spaces are 
defined in terms of comparisons with the model spaces $M_\kappa^2$ (see \cite{BriHae99} 
for more details). Denote $D_\kappa = \pi/\sqrt{\kappa}$. 

\bprop\label{prop-CATk-abconvex}
A CAT$(\kappa)$ space $X$ is $\ds \left(\frac{\mu D_\kappa}{2}, 1 - \cos \frac{\mu \pi}{2}\right)$-convex 
for any $\mu\in(0,1]$.
\eprop
\solution{
Let $x,y,z\in X$ with $d(x,y),d(x,z),d(z,y)\leq \frac{\mu D_\kappa}{2}$ and $t \in [0,1]$. 
By \cite[Lemma 3.3]{Pia11} (see also \cite[Lemma 4.1]{LeuNic13}) we have that
$\ds 
d((1 - t)x + ty,(1-t)x + tz) \le \frac{\sin \frac{t\mu\pi}{2}}{\sin \frac{\mu\pi}{2}} d(y,z).$
Note that
\bua
1-\frac{\sin \frac{t\mu \pi}{2}}{\sin\frac{\mu \pi}{2}} & = & \frac{2\cos\frac{(1+t)\mu\pi}{4}\sin\frac{(1-t)\mu\pi}{4}}{\sin\frac{\mu \pi}{2}}
\ge \frac{2\cos\frac{\mu\pi}{2}\sin\frac{(1-t)\mu\pi}{4}}{\sin\frac{\mu \pi}{2}} \\
& \ge &  (1-t)\cos\frac{\mu\pi}{2}, \quad \text{since } \sin\frac{(1-t)\mu\pi}{4} \geq \frac{1-t}2\sin\frac{\mu\pi}{2}.
\eua
Hence, 
\[\frac{\sin \frac{t\mu\pi}{2}}{\sin \frac{\mu\pi}{2}} \le 1- (1-t)\cos\frac{\mu\pi}{2} =  
t + \left(1 - \cos \frac{\mu \pi}{2}\right)(1-t).
\]
}

We point out that CAT$(\kappa)$ spaces with $\kappa > 0$ do not have in general a convex 
$r$-geodesic bicombing for $r < D_\kappa$ (to see this it suffices to consider the spherical 
space $\mathbb{S}^2$).

Let us recall another notion of convexity for metric spaces, introduced by Ohta \cite{Oht07}.  
Given $L_1,L_2\in[0,\infty)$, a geodesic space $X$ is said to be {\em $L$-convex} for $(L_1,L_2)$ 
if for any $x,y,z \in X$, any constant speed geodesics $\gamma, \xi : [0,1] \to X$ 
with $\gamma(0) = \xi(0) = x$, $\gamma(1)=y$, $\xi(1)=z$ and for every $t\in[0,1]$,
\[
d(\gamma(t), \xi(t))\leq \left(1+L_1\frac{\min\{d(x,y)+d(x,z),2L_2\}}2\right)td(y,z).
\]

An additional related notion says that an $r$-geodesic bicombing on a metric 
space $X$ is {\it weakly convex} if there exists a constant $C \ge 1$ such that
\[d((1-t)x+ty, (1-t)x+tz)\leq Ctd(y,z),\]
for all $t \in [0,1]$ and all $x,y,z\in X$ as in Definition \ref{def-rdelta-Busemann}. 
One can easily see that an $L$-convex space for $(L_1,L_2)$ has a weakly convex $r$-geodesic bicombing, where 
$r\in(0,\infty]$ and $C:=1 + L_1\min\{r, L_2\}$. 

\bfact
We remark that we could have defined an even more general notion: given 
$r>0$ and $\eta:[0,1]\to [0,\infty)$, a metric space with an $r$-geodesic bicombing is {\em $(r,\eta)$-convex} if 
\[
d((1-t)x+ty, (1-t)x+tz)\leq (t+\eta(t))d(y,z),
\]
for all $t \in [0,1]$ and all $x,y,z\in X$ as in Definition \ref{def-rdelta-Busemann}.

This very general definition has the advantage that it covers the case of metric spaces 
with a weakly convex $r$-geodesic bicombing and, thus, of $L$-convex spaces.
\efact

However, we use in this paper Definition \ref{def-rdelta-Busemann}, as this is the notion 
which allows us to get the effective results from the next section.

\section{Effective rates of asymptotic regularity}

Let $r\in(0,\infty], \delta\in[0,1)$, $X$ be an  $(r,\delta)$-convex space,  $C\subseteq X$  a convex subset and 
$T_1, \ldots, T_N:C\to C$ be nonexpansive mappings, where $N\in\Z_+$. 
If $(\lambda_n)$ is a sequence in $[0,1]$ and $u\in C$, one can, obviously, 
define the iterations \eqref{Halpern-N} and  \eqref{alternative-iteration-N} starting with $u$ in this setting, too:
\bea
x_0=u, \quad x_{n+1} = \lambda_{n+1}u + (1-\lambda_{n+1})T_{n+1}x_n,\label{def-Halpern-N-Busemann}\\
x_0=u, \quad x_{n+1} = T_{n+1}(\lambda_{n+1}u + (1-\lambda_{n+1})x_n), \label{def-alternative-N-Busemann}
\eea
where $T_n = T_{n \text{ mod } N}$ and the $\text{mod } N$ function takes values in $1, \ldots ,N$. We  use the following notation $\ds \tnN:= T_{n+N}\cdots T_{n+1}$.

\vspace*{1mm}
The main theorem of the paper is a quantitative result on the asymptotic regularity of the above iterations.

\begin{theorem}\label{main-thm-N}
Let $\eps>0$, $M>0$ be such that $\ds M\leq \frac{r}2$, $\alpha, \gamma:(0,\infty)\to\Z_+$
and $\theta:\Z_+\to \Z_+$. Suppose that
\be
\item $\ds\sum_{n=1}^\infty \lambda_{n}=\infty$  with rate of divergence $\theta$;
\item $\ds \sum_{n=1}^\infty |\lambda_{n+N}-\lambda_n|$  converges with  Cauchy modulus $\gamma$.
\ee
Let 
\bua
\tilde{\Phi}(\eps,M,\gamma,\theta,\delta) & = &  
\theta\bigg(\left\lceil\frac{1}{1-\delta}\right\rceil\left(\gamma\left(\frac\eps{4M}\right)+
\max\left\{\left\lceil\ln\left(\frac{4M}\eps\right)\right\rceil,1\right\}\right)\bigg),\\
\Phi(\eps,M,\gamma,\theta,\delta,N,\alpha) & = & 
\max\left\{\tilde{\Phi}\left(\frac{\eps}2,M,\gamma,\theta,\delta\right), \alpha\left(\frac\eps {4MN}\right)\right\}.
\eua
Assume either
\be
\item\label{main-Halpern-N} $(x_n)$ is given by \eqref{def-Halpern-N-Busemann} with 
$d(x_n,u)\leq M$ for all $n\geq 1$ and $d(u,T_i u)\leq M$ for 
each $i = 1, \ldots, N$, or
\item\label{main-alternative-N}  $(x_n)$ is given by \eqref{def-alternative-N-Busemann} 
with  $d(x_n,u)\leq M$ for all $n\geq 1$.
\ee

Then $\limn d(x_n,x_{n+N})=0$ with rate of convergence $\tilde{\Phi}$. 
Furthermore, if $\limn \lambda_{n+1}=0$ with rate of convergence $\alpha$, then 
$\ds \limn d(x_{n},T_{n,N}(x_n))=0$ with rate of convergence $\Phi$.
\end{theorem}

We give the proof of the theorem in the next section. Let us state now some immediate consequences.

\bcor\label{main-thm-N-bounded}
Let $\eps, M, (\lambda_n), \alpha,\gamma,\theta, \tilde{\Phi},\Phi$ be as above. 
Assume moreover that $C$ is bounded and $M$ is an upper bound on its diameter.

If $(x_n)$ is given by either  \eqref{def-Halpern-N-Busemann} or  \eqref{def-alternative-N-Busemann}, then  
$\ds \limn d(x_n,x_{n+N})=0$ with rate of convergence $\tilde{\Phi}$ and 
$\ds \limn d(x_{n},T_{n,N}(x_n))=0$ with rate of convergence $\Phi$.
\ecor

Thus, for $C$ bounded we obtain a highly uniform rate of asymptotic 
regularity $\Phi$ which does not depend at all on the starting point $u$ and the nonexpansive mappings $T_1,\ldots, T_N$. 
Moreover, the dependence on the set $C$ and the space $X$ is very weak: via $\delta$ and a bound $M \le \frac{r}2$ on the diameter of $C$.

\bcor
Suppose that  $\ds \lambda_n=\frac1{n+1}, \, n\geq 1$. Then $\limn d(x_n,x_{n+N})=\limn d(x_n,T_{n,N}(x_n))=0$ 
with a common rate of convergence 
\[
\Psi(\varepsilon,M,N,\delta) = \exp\left(\left\lceil\frac1{1-\delta}\right\rceil
\left(\left\lceil\frac{8M(N+1)}\eps\right\rceil+2\right)\ln 4\right).
\]
\ecor
\begin{proof}
We can take $\ds \theta(n)=\exp(n\ln4)$, 
$\ds \gamma(\eps)=\left\lceil\frac{N}\eps\right\rceil$ and 
$\ds \alpha(\eps)=\left\lceil\frac{1}\eps\right\rceil$.
\end{proof}

As mentioned before, in the case $N=1$, the iterative scheme defined by 
\eqref{def-Halpern-N-Busemann} yields the usual Halpern iteration, for which rates of 
asymptotic regularity have already been computed in the setting of CAT$(\kappa)$ spaces  
\cite{LeuNic13}. Our main theorem recovers (with a slightly 
modified rate) \cite[Proposition 3.2]{LeuNic13}  since, 
for $\ds M < \frac{D_\kappa}2$, one takes $\ds \mu= \frac{2M}{D_k}$ in Proposition 
\ref{prop-CATk-abconvex} to get that any CAT$(\kappa)$ space is 
$(M, 1-\cos(M\sqrt{\kappa}))$-convex and then apply Corollary \ref{main-thm-N-bounded}. Furthermore, 
we generalize with basically the same bounds the results obtained for the Halpern 
iteration in normed spaces \cite{Leu07a} and, more general, $W$-hyperbolic spaces \cite{Leu09-Habil}.

As we have already pointed out, in this paper we obtain for the first time, even for 
Banach spaces, effective bounds on the asymptotic regularity of the iterations 
\eqref{def-Halpern-N-Busemann} and \eqref{def-alternative-N-Busemann}.  
Recently, using proof mining methods as well, Khan and Kohlenbach \cite{KhaKoh13} obtained in the setting of
uniformly convex Busemann spaces effective results on the asymptotic behavior of a different 
iteration associated to a finite family of nonexpansive mappings which extends the Krasnoselski-Mann iteration of a single 
nonexpansive mapping.
\section{Proof of the main result}

Assume the hypothesis of Theorem \ref{main-thm-N}. As in the case of the Halpern iteration 
associated to a single mapping \cite{Leu07a,Leu09-Habil,KohLeu12a,LeuNic13}, 
we shall apply the following quantitative lemma.

\blem\label{main-lema-as-reg-Halpern-divergent-sum}\cite{LeuNic13,KohLeu12a}
Let $(\alpha_n)_{n\ge 1}$ be a sequence in $[0,1]$ and $(a_n)_{n\geq 1},(b_n)_{n\geq 1}$ be 
sequences in $\R_+$  such that
\beq
a_{n+1}\leq (1-\alpha_{n+1}) a_n + b_n \quad \text{for all~} n\in\Z_+.
\eeq
Assume that $\ds \sum_{n=1}^\infty b_n$  is convergent with Cauchy modulus  $\gamma$ and 
$\ds \sum_{n=1}^\infty \alpha_{n+1}$ diverges with rate of divergence $\theta$.

Then, $\ds\limn a_n=0$ with rate of convergence $\Sigma$ given by 
\beq
\Sigma(\eps,P,\gamma,\theta)=\theta\left(\gamma\left(\frac\eps 2\right)+
\max\left\{\left\lceil\ln\left(\frac{2P}\eps\right)\right\rceil,1\right\}\right)+1
\eeq
where $P>0$ is an upper bound on  $(a_n)$.
\elem

The next lemma is the second main tool for the proof of our main result.

\blem \label{lemma-ineq}
Let $M>0$  satisfy $2M\leq r$. 
\be
\item Assume that $(x_n)$ is given by \eqref{def-Halpern-N-Busemann},  $d(u,T_i u)\leq M$ 
for each $i = 1, \ldots, N$ and $d(x_n,u)\leq M$ for all $n\in\N$. Then, for all $n\ge 1$,
\begin{align*}
d(x_{n+1},T_{n+1}x_{n}) &\le  2M\lambda_{n+1},\\
d(x_n, x_{n+N}) &\le (1-(1-\delta)\lambda_{n}) d(x_{n-1},x_{n+N-1}) + 2M|\lambda_{n + N} - \lambda_n|.
\end{align*}
\item Assume that $(x_n)$ is given by \eqref{def-alternative-N-Busemann} and $d(x_n,u)\leq M$ 
for all $n\in\N$. Then, for all $n\ge 1$,
\begin{align*}
d(x_{n+1},T_{n+1}x_{n}) &\le   M\lambda_{n+1}, \\
d(x_n, x_{n+N}) & \le (1-(1-\delta)\lambda_{n}) d(x_{n-1},x_{n+N-1}) + M|\lambda_{n + N} - \lambda_n|. 
\end{align*}
\ee
\elem
\solution{
(i) First, let us note that $d(u, T_n x_m) \le 2M$ for all $m,n \ge 1$. It follows that
for all $n\ge 1$,
\[
d(x_{n+1},T_{n+1}x_{n}) = \lambda_{n+1} d(u, T_{n+1}x_n) \le 2M\lambda_{n+1} 
\]
and
\begin{align*}
& d(x_n, x_{n+N}) =  d(\lambda_n u + (1-\lambda_n)T_n x_{n-1}, \lambda_{n+N}u + (1-\lambda_{n+N})T_n x_{n+N-1})\\
& \quad \le d(\lambda_{n} u + (1-\lambda_{n})T_n x_{n-1}, \lambda_{n}u+ (1-\lambda_{n})T_n x_{n+N-1})\\
& \qquad + d(\lambda_{n} u + (1-\lambda_{n})T_n x_{n+N-1}, \lambda_{n+N} u + (1-\lambda_{n+N})T_n x_{n+N-1})\\
& \quad \le (1-(1-\delta)\lambda_{n}) d(x_{n-1},x_{n+N-1}) + 2M|\lambda_{n + N} - \lambda_n|.
\end{align*}
(ii) The proof is similar, using that
\begin{align*}
d(x_{n+1},T_{n+1}x_{n}) \le d(\lambda_{n+1}u + (1-\lambda_{n+1})x_n, x_n) \le M\lambda_{n+1} 
\end{align*}
and
\begin{align*}
d(x_{n}, x_{n+N}) & = d(T_n(\lambda_{n} u + (1-\lambda_{n})x_{n-1}), T_n(\lambda_{n+N}u 
+ (1-\lambda_{n+N})x_{n+N-1}))\\
& \le d(\lambda_{n} u + (1-\lambda_{n})x_{n-1}, \lambda_{n+N}u + (1-\lambda_{n+N})x_{n+N-1}).
\end{align*}
for all $n\ge 1$.
}

\subsection{Proof of Theorem \ref{main-thm-N}}

Let $(x_n)$ be given by either  \eqref{def-Halpern-N-Busemann} or  \eqref{def-alternative-N-Busemann}. 
Denote $\ds t_{n+1}= (1-\delta)\lambda_{n} \in [0,1]$. As an immediate consequence of 
Lemma \ref{lemma-ineq}, we get that 
\[d(x_n, x_{n+N}) \le (1-t_{n+1}) d(x_{n-1},x_{n+N-1}) + 2M|\lambda_{n + N} - \lambda_n|.\]
Note that $\ds \sum_{n=1}^\infty 2M|\lambda_{n+N}-\lambda_n|$ converges with Cauchy modulus 
$\ds \tilde{\gamma}(\varepsilon)=\gamma\left(\frac{\varepsilon}{2M}\right)$ and 
$\ds \sum_{n=1}^\infty t_{n+1} = \infty$ with rate of divergence $\ds \tilde{\theta}(n) = 
\theta\left(\left\lceil\frac{1}{1-\delta}\right\rceil n\right)$. 

We apply Lemma \ref{main-lema-as-reg-Halpern-divergent-sum} with $\alpha_n:=t_n$, $P:=2M$, 
$a_n := d(x_{n-1},x_{n+N-1})$ and $b_n := 2M|\lambda_{n+N}-\lambda_n|$ to obtain 
that $\limn d(x_{n-1},x_{n+N-1})=0$ with rate of convergence $\tilde{\Phi}+1$. Hence, 
$\limn d(x_n,x_{n+N})=0$ with rate of convergence $\tilde{\Phi}$.

Assume now that  $\limn \lambda_{n+1}=0$ with rate of convergence $\alpha$. By 
Lemma \ref{lemma-ineq}, it follows that $\ds d(x_{n+1},T_{n+1}(x_{n})) \le 2M \lambda_{n+1}$, hence 

\[ d(x_{n+1},T_{n+1}(x_{n})) \le \frac{\eps}{2N} \text{ for all } n \ge \alpha\left(\frac{\eps}{4MN}\right).\] 

One can easily see that 
\[d(x_{n},\tnN(x_n)) \le  d(x_n, x_{n+N}) + \sum_{i=1}^N d(x_{n+i}, T_{n+i}(x_{n+i-1})).\]

Therefore, $d(x_{n},\tnN(x_n))\le \eps$  for all $\ds n\geq \Phi$. $\hfill \Box$

\mbox{ } 

\noindent
{\bf Acknowledgements:} \\[1mm] 
Lauren\c tiu Leu\c stean was supported by a grant of the Romanian 
National Authority for Scientific Research, CNCS - UEFISCDI, project 
number PN-II-ID-PCE-2011-3-0383. \\[1mm]
Adriana Nicolae was supported by a grant of the Romanian
Ministry of Education, CNCS - UEFISCDI, project number PN-II-RU-PD-2012-3-0152.

\end{document}